\documentclass[12pt, hidelinks]{article}
\usepackage[utf8]{inputenc}
\usepackage{amsmath, amssymb, url, float, graphicx, float}
\usepackage{fullpage}
\usepackage[small,bf]{caption}
\usepackage{cite}
\usepackage{booktabs}
\usepackage{graphicx}
\usepackage{xcolor}

\newcommand{\diam}{\mathbf{diam}}

\newcommand{\E}{\mathbf{E}}

\let\phi\varphi

\newcommand{\ones}{\mathbf 1}
\newcommand{\reals}{{\mbox{\bf R}}}

\DeclareMathOperator*{\argmax}{argmax}

\newcommand{\cf}{{\it cf.}}
\newcommand{\eg}{{\it e.g.}}
\newcommand{\ie}{{\it i.e.}}

\newcommand{\BEAS}{\begin{eqnarray*}}
\newcommand{\EEAS}{\end{eqnarray*}}
\newcommand{\BEA}{\begin{eqnarray}}
\newcommand{\EEA}{\end{eqnarray}}
\newcommand{\BEQ}{\begin{equation}}
\newcommand{\EEQ}{\end{equation}}
\newcommand{\BIT}{\begin{itemize}}
\newcommand{\EIT}{\end{itemize}}

\newif\iftodos

\title{The Specter (and Spectra) of Miner Extractable Value}
\author{
    Guillermo Angeris\thanks{The authors are listed in alphabetical order.}\\\texttt{\small gangeris@baincapital.com}
    \and Tarun Chitra\\\texttt{\small tarun@gauntlet.network}
    \and Theo Diamandis\\\texttt{\small tdiamand@mit.edu}
    \and Kshitij Kulkarni\\\texttt{\small ksk@eecs.berkeley.edu}
}
\date{August 2023}

\begin{document} 
\maketitle 

\begin{abstract}
    Miner extractable value (MEV) refers to any excess value that a transaction
    validator can realize by manipulating the ordering of transactions. In this
    work, we introduce a simple theoretical definition of the `cost of MEV',
    prove some basic properties, and show that the definition is useful via a
    number of examples. In a variety of settings, this definition is related to 
    the `smoothness' of a function over the symmetric group. From this 
    definition and some basic observations, we recover a number of results from 
    the literature.
\end{abstract}

\section*{Introduction}
Decentralized systems, such as public blockchains, allow users to submit
transactions that, if valid, modify the shared state of the system. Examples of 
these transactions include peer-to-peer transfers, liquidations of 
undercollateralized loans, and trades on decentralized exchanges. These
submitted transactions are then aggregated into ordered lists, called
\emph{blocks}, by agents called \emph{validators}. These validators, historically
called \emph{miners}, propose blocks to be included in the blockchain and
generally have the freedom to arrange transactions within a block as they
wish. The validators are privileged agents that can extract a large amount of value
by rearranging transactions in a way that is beneficial to them. For example,
a validator may place their order to purchase some asset before another user's 
purchase order; this action is known as `front-running'. The type of value
that can be extracted by a validator in this way is termed \emph{miner
extractable value} or \emph{MEV}, for short~\cite{daian2020flash}. One can view
this extracted value as excess returns earned by validators within these
systems at the expense of other users, which has amounted to over a billion of
dollars since 2020~\cite{flashbots_mev_explore}. We explore this process in
more detail below.

%

\paragraph{Block construction.} 
Multiple agents are involved in the creation of a block. 
Users send transactions, such as trades from one asset to another on a
decentralized exchange, via a peer-to-peer network. These transactions are received
by validators: network participants who maintain the validity of state in a
public blockchain and
earn fees for producing valid blocks. In most consensus protocols, a validator
has a temporary monopoly over both what transactions are included in a block
and the order in which those transactions are executed. With this additional
power, validators can increase their revenue by reordering or adding 
transactions. For example, a validator may front-run a user's trade, causing 
that user to pay a higher price than they otherwise would.

%

\paragraph{Miner extractable value.} Research into MEV has generally focused on
empirical properties of MEV usually found in practice. The study of MEV, and
its general high-level definition, began empirically
with~\cite{daian2020flash}, which identified ways in which validators could
earn a profit by reordering transactions in ways beneficial to
them. Subsequent work analyzed the profitability of MEV
for particular types of applications~\cite{zhou2021high, qin2022quantifying}.
These analyses focused on identifying particular types of MEV via heuristics and 
estimating the total value `extracted'.

\paragraph{Formalization of MEV.} Some work has formalized the amount of
MEV from reordering transactions in specific settings or for specific
protocols~\cite{kulkarni2022towards,bartoletti2022maximizing}. These works
attempt to formulate optimization problems that characterize MEV in a system, and then prove
bounds on the values returned as solutions. However, these results only apply to 
constant function market makers~\cite{angeris2022}. 
Other forms of extractable value have not been formally
analyzed. For example, no theoretical work formalizes MEV from front-running in 
non-fungible token (NFT) auctions or from liquidations in systems with leverage,
although both of these scenarios generate substantial MEV in practice. 

\paragraph{Mitigation of MEV.} In contrast, many works propose mechanisms to mitigate MEV in the 
general consensus setting~\cite{kelkar2021themis, malkhi2022maximal} or in 
specific applications~\cite{mcmenamin2022fairtradex, greeneoracle, 
chitra2022improving, bebel2022ferveo, johnson2023concave, rondelet2023threshold,
resnick2023contingent, macpherson2023adversarial, xavier2023credible}. These 
mechanisms generally require consensus to enforce extra ordering constraints, 
such as those that respect first-come-first-serve dynamics. To characterize MEV, 
this line of work either considers a specific application, or assumes that MEV
can be reasonably solved by ensuring first-come-first-serve dynamics.

\paragraph{This paper.} In this paper, we formally define the `cost of MEV'
which quantifies the excess amount that a profit maximizing validator can earn
from reordering a set of transactions, relative to a random ordering. We
quantify this profit via a payoff (or value) function that maps a specific
transaction ordering to a profit (or utility) realized by the validator. The
cost of MEV, then, essentially compares profits in a system that is `order
agnostic' versus one that explicitly gives the right to reorder transactions to
a validator who seeks to maximize their payoff. We show that the cost of MEV
shares certain smoothness properties of the underlying value function and show
applications of this definition and consequences via basic examples. We also
connect the notion of smoothness of a value function to the smoothness of
functions defined on graphs, and in particular, the graph's spectral
properties. Our results allow protocol designers and theoreticians to reason
about the cost of MEV in purely abstract terms. Along the way, we derive (and,
in some cases, strengthen) some results already known in the literature as
special cases of the provided definitions and their implications. Though our
definition suggests some requirements for mechanisms and/or applications to
have `low MEV', and we note some of these requirements in this paper, we do not
focus on this line of research, leaving it for future work.

\section{The cost of MEV}
In this section, we define the cost of MEV (miner extractable value) and show
some basic consequences of the definition, including simple bounds on the cost
and the family of functions which saturate such bounds. We then
define a number of special classes of payoffs which satisfy certain properties
with respect to the cost and show basic generalizations of the cost of MEV that
are reasonable.

\paragraph{Set up.} In our set up, there is an action space, defined as a set
$A$, and a value function $f: A^n \to \reals$, where $n$ is the total number of
transactions in the block. The set $A$ may contain actions available to a user,
including making a trade in a decentralized exchange or liquidating an
underwater loan in an on-chain lending protocol. The function $f$ characterizes
the net utility accrued to a transaction validator from a particular ordered
list of transactions. Some examples we will consider later are the payoff of a
sandwich attack on the traders of a constant function market maker and simple
liquidations. We make almost no assumptions about $f$; indeed, the
value function $f$ may be very complicated or difficult to write down in some
practical cases. For concrete examples of such value functions, we encourage
the reader to (temporarily) skip ahead to~\S\ref{sec:smooth} and skim the
`value function' paragraph headings.

\paragraph{Cost of MEV.} We define the \emph{cost of MEV} for a value function $f$ and
for some set of transactions
$x \in A^n$, as
\begin{equation}\label{eq:main}
    C(f, x) = \max_{\pi \in S} f(\pi(x)) - \E_{\pi \sim S} f(\pi(x)).
\end{equation}
Here, $S$ is the set of all possible permutations of length $n$ (also known as
the \emph{symmetric group}) and $\pi \sim S$ means that $\pi$ is uniformly
randomly drawn from $S$. The function $C(f,x)$ characterizes how much the value
of the set of transactions differs between the `worst case' ordering (from the
users' perspective) and the `average case' ordering. If this difference is
large, then a validator may be able to reorder the transactions to extract a
large amount of value from them, relative to an average ordering.
Alternatively, we may view $C(f, x)$ as, roughly speaking, how much an agent
with payoff $f$ would be willing to pay in order to achieve their desired
ordering, versus simply posting the transaction (assuming transactions are
approximately uniformly randomly ordered in arrival times). Note that even
computing this quantity might be difficult since the size of the possible
permutations, $|S|$, is $n!$, which is large even for moderately small $n$, but we
suspect many basic heuristics work very well in practice and assume this
maximum can (at least approximately) be achieved.

For convenience, we will write $\E_\pi$ for $\E_{\pi \sim S}$ and
$\max_\pi$ for $\max_{\pi \in S}$ to avoid overly burdensome notation, when the
meaning is clear or easily implied. 

\paragraph{Homogeneity.} Note that the cost of MEV is homogeneous
in the value function in that, for any $\alpha \ge 0$, we have
\begin{equation}\label{eq:homogen}
    C(\alpha f, x) = \alpha C(f, x).
\end{equation}
Because of homogeneity, it often suffices to pick functions $f$ that are `bounded' or
normalized in some sense, for fixed $x$ (\eg, $\max_\pi |f(\pi(x))| = 1$), since we can
always rescale the result as needed.

\paragraph{Translation invariance.} The function $f$ and any `translation' of this
function, $\tilde f(x) = f(x) + \alpha$ both have the same cost of MEV,
\[
    C(\tilde f, x) = C(f, x),
\]
for any $\alpha \in \reals$.

\paragraph{Normalization.} Translation invariance, together with homogeneity,
means that it usually suffices to consider \emph{normalized} functions that
satisfy $0 \le f \le 1$ with $\max_\pi f(\pi(x)) = 1$, to derive most general
properties of the cost of MEV. Since any function can be scaled and translated
to satisfy these bounds, we may set
\[
    \tilde f(x) = \frac{f(x) - \min_\pi f(\pi(x))}{\max_\pi f(\pi(x)) - \min_\pi f(\pi(x))},
\]
if $f(\pi(x))$ is not constant over $\pi$, for given $x$, and otherwise $\tilde
f(x) = 1$. Clearly $0 \le \tilde f(x) \le 1$ and we have
\[
    C(f, x) = (\max_\pi f(\pi(x)) - \min_\pi f(\pi(x))) C(\tilde f, x).
\]

\paragraph{Subadditivity.} The cost of MEV is subadditive with respect to the
value functions. In particular, given two value functions $f, g: A^n \to
\reals$ we have
\[
    C(f+g, x) \le C(f, x) + C(g, x),
\]
for any transaction $x \in A^n$, which follows from the fact that the $\max$ is
subadditive and the expectation is linear. Combined with the fact that $C$ is
homogeneous in its first argument, subadditivity implies that $C$ is sublinear with
respect to its first argument.

\paragraph{Permutation invariance.} We also have the property that,
for any permutation $\pi \in S$, the cost of MEV is unchanged; \ie,
\[
    C(f\circ \pi, x) = C(f, x),
\]
and
\begin{equation}\label{eq:perm}
    C(f, \pi(x)) = C(f, x).
\end{equation}
This happens because $S$ is a group, so $\pi S = S$ for any $\pi \in S$. In
this way, it makes sense to consider the set of possible transactions $A^n$
modulo the set of possible permutations $S$, denoted $A^n/S$ as the `natural
set' over which the cost $C$ is measured, as opposed to measuring the cost directly over
transactions $x \in A^n$. We do not take this approach here, but do use these
properties later in this section to characterize the worst-cost functions.

\subsection{Worst-cost functions}
In this section we introduce the worst-cost functions: a set of payoff
functions that are similar in spirit to the payoffs that liquidators achieve on
many decentralized finance protocols. (In some special cases, worst-cost
functions can be written as a type of `liquidation-like' payoff, as in the
example provided.) We show that such functions are exactly those that saturate
a global bound on the cost of MEV and that they form a basis: every payoff
function $f$ can be written as a linear combination of worst-cost functions.

\subsubsection{Basic bounds}
For a fixed set of transactions $x$, it is not hard to show that, over the set
of normalized functions $0 \le f(\pi(x)) \le 1$ for all $\pi \in S$, where there
exists a $\pi$ such that $f(\pi(x)) = 1$, we have that
\begin{equation}\label{eq:bound}
    C(f, x) \le 1 - \frac{|F(x)|}{n!},
\end{equation}
here $F(x)$ denotes the fixed points of $x$ under the action of $S$,
\[
    F(x) = \{\pi \in S \mid \pi(x) = x\}.
\]
(Such fixed points are sometimes known as the \emph{stabilizers} of $x$ by
$S$.) To see bound~\eqref{eq:bound}, note that the maximum over $\pi \in S$ of
$f(\pi(x))$ is achieved at some $\pi$. But, for any $\pi' \in F(x)$ we have
$\pi(\pi'(x)) = \pi(x)$, so the maximum, which is at most $1$ using our
normalization, is achieved at least $|F(x)|$ permutations. This means, since
$f(\pi(x)) \ge 0$ for all $\pi \in S$,
\[
    \E_\pi[f(\pi(x))] \ge \frac{|F(x)|}{n!},
\]
and the bound~\eqref{eq:bound} follows. Note that the set $F(x)$ is nonempty
since the trivial permutation is always included, so $|F(x)| \ge 1$.

\paragraph{Orbits.} A deep---if easy to prove---statement from group theory,
sometimes called the orbit-stabilizer theorem~\cite[\S2]{clark1984elements}, is
that
\[
    \frac{|F(x)|}{n!} = \frac{1}{|S(x)|},
\]
where $S(x) = \{\pi(x) \mid \pi \in S\}$ is the image of $x$ under all possible
permutations $S$. (This set $S(x)$ is also known as the \emph{orbit} of $x$
under the action of the symmetric group.) We then have the bound
\[
    C(f, x) \le 1 - \frac{1}{|S(x)|}.
\]
In a certain sense, we may view $S(x)$ as the size of the permutations that
actually matter, since some permutations simply leave $x$ unchanged, and
each of these permutations is counted only once in the set $S(x)$.

\paragraph{Implications.} Given some $x \in A^n$ with $x_i \ne x_j$ for all
$i\ne j$, any permutation $\pi \in S$, except the identity (which leaves
the order unchanged) will have $\pi(x) \ne x$ so $|S(x)| = n!$. We recover
the `simpler' result
\[
    C(f, x) \le 1 - \frac{1}{n!}.
\]
More generally, if the function $f$ is not normalized, but is nonnegative, we have
the following bound, using~\eqref{eq:bound} and the homogeneity of $C$ in its
first argument~\eqref{eq:homogen},
\[
    C(f, x) \le \max_\pi f(\pi(x)) \left(1 - \frac{|F(x)|}{n!}\right).
\]

\paragraph{Basic example.} We consider a toy example that is a liquidation-like
payoff. Let the action space $A = \{+1, -1, L\}$, where $L$ is a symbol
indicating `liquidate', while $\pm 1$ denote `trades' made with some market
that increase (or decrease) the price by one unit. The payout to the
liquidator, which we assume to be the validator, is
\[
    f(x) = \begin{cases}
        1 & \text{there exists a}~k~\text{such that}~x_1 + \dots + x_k \ge n/2 ~\text{and}~x_{k+1} = L \\
        0 & \text{otherwise}.
    \end{cases}
\]
(In the first case, we implicitly assume the first $k$ trades are
$\pm 1$.) We may interpret this function as: the validator may insert a
liquidation at any point in time; they succeed at liquidating only when the
price of the asset reaches a certain threshold and their liquidation happens
before any other trades lower the price of the asset below their threshold.
Assuming that $n$ is even and the (fixed) set of transactions $x$ is
\[
    x = (1, \dots, 1, L, -1, \dots, -1),
\]
with $n/2$ entries equal to $+1$ and $n/2-1$ entries equal to $-1$, then,
clearly $f(x) = 1$, so $\max_\pi f(\pi(x)) = 1$. On the other hand, computing
$\E_\pi f(\pi(x))$ is slightly trickier. Note that $f(\pi(x)) = 1$ if, and only
if, the permutation leaves the first $n/2$ elements equal to $1$, the
$(n/2+1)$st element equal to $L$, and the remaining $(n/2 - 1)!$ elements equal
to $-1$, in other words, when the permutation $\pi$ leaves $x$ unchanged,
$\pi(x) = x$. This happens a total number of $(n/2)!(n/2-1)!$ ways, which gives
\[
    |F(x)| = (n/2)!(n/2-1)!,
\]
or, since $f(\pi(x))$ is exactly one when $\pi(x) = x$, \ie, when $\pi \in F(x)$,
and zero otherwise, we have
\[
    \E_\pi f(\pi(x)) = \frac{(n/2)!(n/2-1)!}{n!} = \binom{n}{n/2}^{-1}\frac{2}{n},
\]
which saturates the bound~\eqref{eq:bound}. (As we will see next, this is one
example of a `worst-cost function'.) Putting it all together,
\begin{equation}\label{eq:cost-example}
    C(f, x) = 1 - \binom{n}{n/2}^{-1}\frac{2}{n} ~ \sim ~ 1 - 2^{-n},
\end{equation}
which is very close to $1$ for even moderate values of $n$. (The right hand
side follows by using the approximation $n! \sim n^n$, and the result is off by
no more than a polynomial factor of $n$.)

This example is easily generalized to any sequence of actions $x \in A^n$ such
that exactly one of the following three conditions is true: either $x_i > 0$,
$x_i < 0$, or $x_i = L$ (for at most one index $i$) and we allow any
liquidation threshold price $p$ (compared to just $n/2$) satisfying
\[
    \sum_{i\mid x_i > 0} x_i + \max_{i\mid x_i < 0} x_i < p \le \sum_{i\mid x_i > 0} x_i,
\]
where we ignore the single term with index $k$ such that $x_k = L$ in the sum.
In this case, the bound is very similar to the previous:
\[
    C(f, x) = 1-\frac{m!(n-m-1)!}{n!} = 1 - \binom{n}{m}^{-1}\frac{1}{n-m},
\]
where $m$ is the number of indices $i$ such that $x_i > 0$.

\subsubsection{Worst-cost functions}
A particularly interesting set of functions to analyze are those which
saturate the bound~\eqref{eq:bound}. We discuss such functions
(and some consequences) here.

\paragraph{Fixed transactions.} If a function achieves equality for the
bound~\eqref{eq:bound}, for fixed $x \in A^n$, we will call it a
\emph{worst-cost function for fixed transactions} $x$. The maximum here is
achieved exactly by functions of the form
\begin{equation}\label{eq:ind}
    f(\pi(x)) = \ones[\pi(x) = y],
\end{equation}
where $y \in S(x)$ is some fixed element and $\ones[\pi(x) = y]$ is the 0-1
indicator function that is 1 if $\pi(x) = y$ and 0, otherwise. To see that any
(normalized) nonzero function that meets the bound must be of this form, note
that $f(\pi^\star(x)) = 1$ at some $\pi^\star \in S$, so
\[
    f(\pi^\star\circ\pi'(x)) = f(\pi^\star(x)) = 1,
\]
for every $\pi' \in F(x)$, where the first equality follows since $F(x)$ are
exactly the permutations $\pi'$ for which $\pi'(x) = x$. This implies that
\[
    \E_\pi [f(\pi(x))] \ge \frac{|F(x)|}{n!},
\]
with equality only when $f(\pi(x))$ is zero at all permutations for which
$\pi(x)\ne \pi^\star(x)$. Now, if $f$ meets the bound~\eqref{eq:bound} at
equality, then
\[
    C(f, x) = 1 - \frac{|F(x)|}{n!},
\]
but, since $f$ is normalized then $\max_\pi f(\pi(x)) = 1$ so
\[
    \E_\pi[f(\pi(x))] = \frac{|F(x)|}{n!},
\]
which, from before, implies that $f(\pi(x)) = 0$ except when $\pi(x) =
\pi^\star(x)$, where it is 1; \ie, that $f$ can be written as
\[
    f(\pi(x)) = \ones[\pi(x) = \pi^\star(x)],
\]
setting $y = \pi^\star(x) \in S(x)$ completes the claim.

\paragraph{General worst-cost functions.} Showing that worst-cost functions exist for
any $x$, rather than fixing $x$ ahead of time, is similar, but slightly more
involved than the previous case. If $0 \le f \le 1$, then the bound
of~\eqref{eq:bound} still holds, but defining the class of functions $f$ that
meet this bound for all possible transactions $x$ is somewhat more delicate. 

In this case, we define the set of equivalence classes of the action space
$A^n$ under the possible permutations, which we denote
\[
    A^n / S = \{S(x) \subseteq A^n \mid x \in A^n \},
\]
Note that the sets $A^n/S$ form a partition of $A^n$ as they are pairwise
disjoint and their union is all of $A^n$, which follows from the fact that $S$
is a group. Additionally, since $x \in S(x)$, each set is nonempty. The set of
functions that saturate the bound for all $x \in A^n$, \ie, the global
\emph{worst-cost functions} are then exactly those of the form
\begin{equation}\label{eq:worst-cost}
    f(x) = \sum_{Q \in A^n/S} \ones[x = x_Q],
\end{equation}
where $x_Q \in Q$ for each $Q \in A^n/S$. This (potentially uncountable) sum is
justified since, for any $x$, there is at most one term that will be nonzero.
We can view the function $f$ as picking one `canonical' permutation for every
equivalence class $Q \in A^n/S$ and assigning it a value of $1$ if $x \in Q$
and matches this permutation exactly. Otherwise, it assigns a value of 0. That
this family of functions is nonempty, for general sets $A^n$, requires the
axiom of choice, but these sets are generally far more structured in practice,
so it should be possible to construct this worst-cost function directly.

\paragraph{Proof.} The proof that every normalized worst-cost function is
exactly of this form follows from viewing $f$ restricted a single equivalence
class $Q \in A^n/S$. Let $f$ be a function of the form
of~\eqref{eq:worst-cost}, then, for every $x \in Q$,
\[
    f(x) = \ones[x = x_Q],
\]
since $x$ is in exactly one equivalence class $Q$. Of course, by
definition of $Q$, we have that $\pi(x) \in S(x) = Q$ for every $\pi \in S$, so
\[
    f(\pi(x)) = \ones[\pi(x) = x_Q],
\]
which we know are exactly the worst-cost functions for a fixed $x \in A^n$ from
the previous discussion. Since
\[
    C(f, x) = C(f, y)
\]
for any $y \in S(x) = Q$ by using the permutation invariance of
$C$ in~\eqref{eq:perm}, then this is also a worst-cost function for any $y \in Q$.
Finally, summing over all possible $Q$ and noting that at most one such
indicator has value 1 for each $Q$, and therefore for each $x \in A^n$, gives the
final result.

\paragraph{Worst-cost functions as a basis.} An interesting fact is that the
set of all worst-cost functions is a basis for all functions $f$. In
particular, we have that $f(\pi(x))$, viewed as a function of $\pi$ and holding
$x \in A^n$ constant, can be written as a linear combination of worst-cost
functions:
\begin{equation}\label{eq:basis}
    f(\pi(x)) = \sum_{y \in S(x)} f(y) \ones[\pi(x) = y].
\end{equation}
This expression is well-defined since $S(\pi(x)) = S(x)$ for any permutation
$\pi$, by definition. It is not hard to see that equation~\eqref{eq:basis} is
true. Consider $f(\pi(x))$, then $\pi(x) \in S(x)$ so exactly one term in the
sum is nonzero: the one corresponding to $f(\pi(x))$. The argument is easy to
generalize to the broader case where the transactions $x$ are not fixed and in
that case we have
\[
    f(x) = \sum_{Q \in A^n/S} \sum_{y \in Q} f(y)\ones[x = y],
\]
and we can see that this indeed corresponds to $f$ by (carefully!) interpreting
the sums.

\subsection{Fair and unfair functions}
In this section we define perfectly fair functions, which is a strong
requirement that a few value functions meet in practice. We show that unfair
functions---those which have high cost of MEV---satisfy a number of bounds of
possible interest. For example, if a normalized function is very `localized'
(\ie, it takes on very large values at only a small number of positions and
small values elsewhere) then it is very unfair. We also show a partial
converse: any `large enough' function that is has high cost must be reasonably
localized.

\paragraph{Perfectly fair functions.} Since the maximum is always no smaller
than the expectation, we have that
\[
    C(f, x) \ge 0.
\]
We say a value function $f$ is \emph{perfectly fair for transactions} $x\in A^n$ if
\[
    C(f, x) = 0.
\]
Note that $f$ is perfectly fair for transactions $x$ if, and only if,
\[
    f(\pi(x)) = f(x),
\]
for all $\pi \in S$. Based on this, we will say $f$ is \emph{perfectly fair} if
\[
    C_s(f) = \sup_{x \in A^n} C(f, x) = 0,
\]
since $C_s(f) = 0$ implies that $C(f, x) = 0$ for all $x \in A^n$.
Equivalently, the perfectly fair functions are the set of symmetric functions,
defined as those which have $f \circ \pi = f$ for all $\pi \in S$. While this
might seem like an overly-burdensome restriction, there do exist examples of
perfectly fair functions in practice, these include those of concave pro-rata
games~\cite{johnson2023concave} or sealed-bid auctions~\cite{roughgarden2016},
both of which are mechanisms that do not depend on the order in which user
actions were received. We will see how to construct a `perfectly fair' function
(in expectation) from a general payoff function~$f$ in the extensions.

\paragraph{Unfair functions.} We can push some of the techniques described here
to create lower bounds for general functions $f$ and fixed transactions $x$. In
particular, if a bounded function $f \le 1$ satisfies $f(\pi(x)) \ge \beta$ for
some $\pi \in H \subseteq S$ and $f(\pi(x)) \le \alpha$ for each $\pi \in S
\setminus H$ with $0 \le \alpha \le \beta \le 1$, then $\max_\pi f(\pi(x)) \ge
\beta$ while
\[
    f(\pi(x)) \le \ones[\pi \in H] + \alpha \ones[\pi \in S \setminus H],
\]
for every $\pi \in S$, so
\[
    \E_\pi[f(\pi(x))] \le \alpha + \frac{|H|}{n!}(1-\alpha).
\]
Putting this all together gives the following lower bound on $C(f, x)$,
\begin{equation}\label{eq:set-bound}
    C(f, x) \ge \beta - \alpha - \frac{|H|}{n!}(1-\alpha).
\end{equation}
We can think of the set $H\subseteq S$ as the set of permutations under which
$f$ is `spiky' for transactions $x$. If this set $H$ is (a) small relative to
$S$ and (b) has values that are bounded away from the remaining permutations
with a large margin $\beta \gg \alpha$, for instance, then the cost of MEV is
always guaranteed to be large. 

A particularly interesting case to analyze is when $\alpha = 0$ and $f(\pi(x))$
is nonconstant over $\pi$. In this case, $f(\pi(x))$, seen as a function over
$\pi \in S$, has support included in $H$ and, since the maximum is achieved
somewhere in $H$, we must have that $\beta = 1$, assuming $f$ is normalized.
This gives the simple bound:
\[
    C(f, x) \ge 1 - \frac{|H|}{n!},
\]
Note that this bound is tight against~\eqref{eq:bound} exactly when $H = F(x)$,
the set of stabilizers of $x$ under $S$, which is, when $f$ is normalized,
whenever $f(\pi(x)) = \ones[\pi(x) = \pi'(x)]$ for some fixed permutation $\pi'
\in S$, giving another proof for the form of worst-cost functions. So, if the
support is smaller than the set of all permutations, $|H| < n!$, then
necessarily we have positive cost of MEV. If the function $f$ is not normalized
but is nonnegative with $f(\pi(x)) \ge 0$ for all $\pi \in S$, we have that
\[
    C(f, x) \ge \left(1 - \frac{|H|}{n!}\right)\max_{\pi \in S} f(\pi(x)),
\]
by the homogeneity of $C$ in its first argument~\eqref{eq:homogen}.


\paragraph{Partial converse.} We can also provide a partial converse to the
above that says: if a normalized function has high cost of MEV, then the function
`almost looks like' a worst-cost function. More explicitly, we have the
following two claims. First, if $C(f, x) \ge \alpha$, for a fixed list of
transactions $x$, then
\[
    f(\pi(x)) \ge \alpha\ones[\pi(x) = y],
\]
for some $y \in S(x)$; that is, the function $f$ is at least as large as some
indicator. This result follows pretty much directly from the definition: there
is some $\pi$ such that $f(\pi(x)) \ge \alpha$ and setting $y = \pi(x) \in
S(x)$ suffices. Of course, the function $f$ could potentially be no smaller
than an indicator over a much larger set, call it $T \subseteq S(x)$, but the
following statement bounds the size of the resulting set. If, for some $\eta >
0$ we have that
\[
    f(\pi(x)) \ge \eta \ones[\pi(x) \in T],
\]
then the size of the set $T$ must be bounded by
\[
    |T| \le \frac{(1-\alpha)n!}{\eta |F(x)|}.
\]
That is, fixing $\eta > 0$: if the cost, bounded by $\alpha$, is very
large ($\alpha \approx 1$) then the set $T$ must be very small for reasonable
values of $\eta$. Another way of stating this is: if the function has high cost,
then its largest values are very localized. (Where the most localized
functions with $|T| = 1$ are exactly the worst-cost functions.) To see this, if
$f$ is normalized and $f(\pi(x)) \ge \eta \ones[\pi(x) \in T]$, we have
\[
    \alpha \le C(f, x) \le 1 - \frac{1}{n!}\eta |F(x)||T|,
\]
and the result follows by rearranging.

\subsection{Discussion and generalizations}

From the previous discussion on the cost of MEV, we may view the cost of MEV,
$C$, as measure of how `spiky' or `localized' a particular payoff function $f$
is, with respect to the permutations; the higher the cost, the more localized
such a function is, and vice versa. (We will see another connection to this
interpretation, via the Fourier transform on graphs, in a later section.) We
may also view the cost of MEV as answering the question: how much is a user
willing to pay (in expectation) in order to ensure that their preferred
ordering is achieved on chain? Alternatively, we may view the function $C$ as a
type of `efficiency' metric for MEV, when the functions are normalized: values
very close to 1 denote a near-maximal amount of miner extractable value
possible and values very close to 0 denote small amounts. Of course, there are
many other reasonable measures of the cost of MEV that are very similar in
spirit to the one presented above. We outline some here.

\paragraph{Ratio.} Another reasonable metric to measure the cost of MEV is
a `ratio' as opposed to an absolute difference:
\[
    \tilde C(f, x) = \frac{\E_\pi f(\pi(x))}{\max_\pi f(\pi(x))}.
\]
For normalized functions this reduces to $\E_\pi f(\pi(x))$ and is homogeneous
of degree zero in its function argument. Note also that it is at most 1 for any
function $f$, and is bounded between 0 and 1 for nonnegative $f$. Many results
have direct analogues in this case, but we find that working with the additive
form is simpler and has a number of other interesting results. For example,
this multiplicative variant, $\tilde C$, is not translation invariant, nor does
it satisfy the smoothness results presented in the next section.

\paragraph{Randomness.}
Many papers have pointed out the importance of randomness in achieving
potentially fair mechanisms. This requires extending the current definition
slightly to functions $f: A^n \times \Omega \to \reals$ where $\mathbf{P}$ is a
probability distribution over the sample space $\Omega$, in the following way:
\[
    \bar C(f, x) = \max_{\pi \in P} \E_{\omega \sim \mathbf{P}} f(\pi(x),
    \omega) - \E_\pi[\E_{\omega\sim \mathbf{P}}f(\pi(x), \omega)].
\]
The order of the operations has the following interpretation: the mechanism
uses some randomness, unknown to the user until the ordering $\pi$ has been
provided, and the user wishes to compare the expected value of an ordering that
maximizes their cost function when compared to a uniformly chosen ordering.
This is just a special case of~\eqref{eq:main} where the function being
evaluated is $\bar f: A^n \to \reals$ where $\bar f$ is
\[
    \bar f(x) = \E_{\omega \sim \mathbf{P}} f(x, \omega),
\]
though the practical interpretation differs substantially.

\paragraph{An aside: constructing `fair' functions.} With this extension, we can
convert any (deterministic) function $f: A^n \to \reals$ to a `perfectly fair'
function, in expectation, by letting the sample space be the set of
permutations, $\Omega = S$, letting the distribution be uniform over the sample
space, $\mathbf{P}(\omega) = 1/n!$, and by setting
\[
    \tilde f(x, \omega) = f(\omega(x)).
\]
From before, this function has the interpretation that the protocol chooses a
uniformly random permutation $\omega$ after any set of transactions $x$ has
been provided and permutes these transactions with the randomly chosen
permutation to get $\omega(x)$. To prove fairness, note that, for any $\pi \in
P$,
\[
    \E_{\omega \sim \Omega} f(\omega(\pi(x))) = \E_{\omega \sim \Omega}f(\omega(x)),
\]
since $\omega \circ \pi$ is also uniform over $S$, so
\[
    \bar C(f, x) = 0.
\]
and therefore
\[
    \bar C_s(f) = \sup_{x \in A^n} \bar C(f, x) = 0.
\]
for any $x \in A^n$. 


\paragraph{Local vs.\ global bounds.} In a certain sense, lower bounds on $C_s$
are easy: show that there exist some set of transactions for which the cost $C$
is large. We may then view lower bounds as `local' in that a single example
suffices to show that a lower bound on the cost holds. Upper bounds are,
on the other hand, much harder, as we must make a claim about all possible
lists of transactions from the set of actions $A^n$, which we may view as a
type of global statement.

\paragraph{Going forward.} From here, it is not clear that a lot more can be
said in a general setting without one of two assumptions: either (a) assume
some structure about the set $A$ (such as, \eg, metric structure, or some
notion of smoothness, which will play in various ways with the symmetric
group), or (b) assume that the transactions $x \in A^n$ are fixed and, instead,
look at the properties of $C$ from a `local' perspective. This latter view has
deep connections to spectral graph theory and Fourier transforms over groups.
In a sense, we can view the former statement as being what we can say over all
transactions, whereas the latter property is what we can say over all
permutations. We describe each approach in the following two sections.

\section{Smoothness and permutations}
We now introduce the `global' approach to analyzing the cost of MEV of a value
function over all transactions in $A^n$. To do this, we use smoothness
properties of the value function and some metric structure on the set of
transactions. We will assume that there is a subset $B \subseteq A^n$ closed
under permutations, \ie, $\pi B = B$ for all $\pi \in S$ and that $f: B \to
\reals$. If $B$ can be endowed with a metric $d: B \times B \to \reals_+$ that
is \emph{permutation independent}, that is, if, for any $\pi \in S$, we have
\[
    d(\pi(x), \pi(y)) = d(x, y),
\]
then there are useful smoothness results that hold for $C$. (Examples of
permutation-independent metrics include those induced by the standard norms,
for example.)

\subsection{Smoothness bounds}
In general, $C$ has similar smoothness characteristics
as the function $f$ over $B$. For example, if $f$ is $L$-Lipschitz-continuous
over $B$, that is, if for $L\ge 0$, we have
\[
    |f(x) - f(y)| \le L d(x, y),
\]
for any $x, y \in B$, then $C(f, \cdot)$ is $2L$-Lipschitz-continuous in that,
\begin{equation}\label{eq:smooth-bound}
    |C(f, y) - C(f, x)| \le 2Ld(x, y).
\end{equation}

\paragraph{Proof.} This is easy to see since, for the expectation, we have
\[
    |\E_{\pi'} f(\pi'(x)) - \E_\pi f(\pi(y))| \le \E_\pi |f(\pi(x))-f(\pi(y))| \le \E_\pi[Ld(\pi(x), \pi(y))] = Ld(x, y),
\]
while for $\max$ we have
\[
    f(\pi'(x)) - \max_{\pi} f(\pi(y)) \le \max_\pi |f(\pi(x)) - f(\pi(y))| \le L\max_{\pi}d(\pi(x), \pi(y)) = Ld(x, y)
\]
for all $\pi' \in S$ so
\[
    |\max_\pi f(\pi(y)) - \max_{\pi'} f(\pi'(x))| \le Ld(x, y).
\]
Applying the triangle inequality and using the definition of $C$ yields the
final result~\eqref{eq:smooth-bound}.

\paragraph{Tightness.} This bound is tight in that the constant cannot be
improved for general functions $f$, sets $B$, and metrics $d$, unless other
assumptions are made. We can see this from the simple example where $B = A^n =
\reals^n$, the metric $d$ is the one implied by the $\ell_1$ norm, and the
value function is
\[
    f(x) = x_1 - x_2 - \dots - x_n.
\]
First, we see that $f$ is $1$-Lipschitz continuous in the $\ell_1$ norm,
$\|\cdot\|_1$, that is:
\[
    |f(y) - f(y) | \leq 1 \|x- y\|_1 = \sum_{i=1}^n |x_i - y_i|,
\]
since
\[
    |f(x) - f(y)| = |(x_1 - y_1) - (x_2 - y_2) - \dots - (x_n - y_n)| \le \|x - y\|_1,
\]
by applying the triangle inequality $n$ times. Now, consider $x = 0$ and $y =
e_1 = (1, 0, \dots, 0)$. Since, $\|x - y\|_1 = \|y\|_1 = 1$ the smoothness
bound~\eqref{eq:smooth-bound} implies that
\begin{equation}\label{eq:ex-bound}
    |C(f, x) - C(f, y)| \le 2\|x - y\|_1 = 2,
\end{equation}
where we have used the fact that $f$ is 1-Lipschitz continuous with respect to
the $\ell_1$-norm. We will show that this bound on the
cost~\eqref{eq:ex-bound}, implied by the smoothness
bound~\eqref{eq:smooth-bound} is saturated in a certain limit. First,
consider $f(\pi(x))$ where $x = 0$, where we have
\[
    f(\pi(x)) = 0
\]
for any $\pi \in S$ so $C(f, x) = 0$. Now, consider $y = e_1$. Here, where we
have $\max_{\pi} f(\pi(y)) = 1$ and
\[
    \E_{\pi} f(\pi(y)) = \frac{1}{n} - \left(1-\frac1n\right) = \frac{2}{n} - 1.
\]
We can see this since $f(\pi(y)) = 1$ only when when $\pi$ fixes the first
element of $y$ to be $1$, which happens with probability $1/n$, and takes the
value $-1$ at all other times, with probability $1-1/n$. Using the definition
of $C(f, y)$ we then get
\[
    C(f, y) = 2\left(1 - \frac{1}{n}\right),
\]
so
\[
    |C(f, x) - C(f, y)| = |C(f, y)| = 2\left(1 - \frac{1}{n}\right).
\]
(Compare this with the bound~\eqref{eq:ex-bound}.) Taking $n \uparrow \infty$,
we see that the cost of MEV approaches $2$, which saturates the bound provided
by~\eqref{eq:smooth-bound} in the limit. (Scaling this choice of $f$ by any
value $L > 0$ shows that this bound is saturated for any $L$.)

\subsubsection{Global upper bounds}
If we know the function is smooth and the set of transactions $B$ is bounded
with diameter at most $t$, \ie,
\[
    t \ge \sup_{x, y \in B} d(x, y),
\]
then we have the immediate upper bound on the maximum cost
\[
    C_s(f) = \sup_{x \in B} C(f, x),
\]
given by:
\begin{equation}\label{eq:smooth-upper-bound}
    C_s(f) \le 2Lt + \inf_{x \in B} C(f, x).
\end{equation}
So, from the smoothness of $f$ we receive a general bound on the cost of MEV
that matches our intuition: if the function $f$ is `smooth' over the
transactions, even independent of possible ordering effects, then it has `low'
cost. Additionally, if any action $x \in B$ has $x_i = x_j$ for every $i, j =1,
\dots, n$, then the second term in the right-hand-side
of~\eqref{eq:smooth-upper-bound} is zero, so
\[
    C_s(f) \le 2Lt.
\]
The bound~\eqref{eq:smooth-upper-bound} generalizes and improves upon some
known results over constant function market
makers~\cite{chitraDifferentialPrivacyConstant2022, kulkarni2022towards}, and
we show some explicit examples below.

\subsubsection{Extensions}
Note that we do not use the fact that $d$ is a metric
in any of the above derivations, so it may be any function which satisfies the
Lipschitz-like condition
\[
    |f(x) - f(y)| \le L d(x, y).
\]
Additionally if $d$ is not permutation independent, we may define
\[
    \tilde d(x, y) = \max_\pi d(\pi(x), \pi(y)),
\]
which is permutation independent and satisfies $\tilde d \ge d$. All bounds
above hold over this new `distance-like' function, $\tilde d$.

\subsection{Smooth mechanism examples}\label{sec:smooth}
We show some basic applications of the smoothness upper
bound~\eqref{eq:smooth-upper-bound} to frontrunning and sandwiching in
decentralized exchanges.

\subsubsection{Frontrunning}
A classic example of MEV is the notion 
of \emph{frontrunning}; \ie, the fact that parties (such as validators)
are willing to pay some cost in order to be the first transaction executed.
\paragraph{Action space.} In this particular case, we will assume that the set
of possible actions is
\[
    A = \underbrace{(\reals_+ \times \{0\})}_\text{actions by traders} \cup \underbrace{(\{0\} \times \{\delta\})}_\text{actions by validator},
\]
where $\delta \in \reals_+$ is some nonnegative number. We can view a tuple
$(x, y) \in A$ as an action, where $x$ (a trade made by a trader) is nonzero
only when $y$ (a trade made by the validator) is zero, and $y = \delta$ is the
only possible trade that can be made by a validator. We then constrain the
possible actions to be bounded in the following sense:
\[
    B = \{((x_1, y_1), \dots, (x_n, y_n)) \in A^n \mid \|x\|_1 \le M ~ \text{and exactly one entry of $y$ is nonzero}\}.
\]
In other words, the total volume traded is constrained to be at most $M$,
while the validator only submits exactly one trade of size $\delta$. To avoid
overly burdensome notation, we will overload symbols slightly and write
$(\tilde x, \tilde y) \in B$ to mean
\[
    ((\tilde x_1, \tilde y_1), \dots, (\tilde x_n, \tilde y_n)) \in B.
\]
In other words, $\tilde x$ is the vector of trades by traders and $\tilde y$ is the
vector of trades by the validator. By definition, $\tilde x_i$ is nonzero
only when $\tilde y_i = 0$ and vice versa.

\paragraph{Metric.} In this case, we will simply set the `metric' $d$ to be
given by the following function, defined for some
\[
    z = (\tilde x, \tilde y) \in B,
\]
as
\[
    d(z, z') = \sum_{i=1}^n \max\{|\tilde x_i|, |\tilde x_i'|, |\tilde x_i - \tilde x_i'|\},
\]
where $\tilde x'$ has the same meaning as $\tilde x$, except for action $z'$.
This `metric' is permutation-independent and is not a true metric. It measures
something very similar, but not quite equivalent, to the total volume of trades
on the CFMM made by traders.

\paragraph{Diameter bound.} Note that $B$, the set of allowable actions, has
`diameter' bounded with respect to this `metric' using the definition of the
set $B$ and so satisfies 
\begin{equation}\label{eq:diam-bound}
    d(z, z') \le 2M
\end{equation}
for every $z, z' \in B$. To see this, apply the triangle inequality,
\[
    |\tilde x_i - \tilde x_i'| \le |\tilde x_i| + |\tilde x_i'|,
\]
which implies that
\[
    d(z, z') \le \sum_{i=1}^n \max\{|\tilde x_i|, |\tilde x_i'|, |\tilde x_i| + |\tilde x_i'|\} = \sum_{i=1}^n |\tilde x_i| + |\tilde x_i'| \le 2M.
\]
Surprisingly, this bound is tight when the number of trades is greater than
one. We can see this by setting $\tilde x_i = 2M/n$ for even indices $i$ and
$0$ otherwise, and $\tilde x_i' = 2M/n$ for odd indices $i$ and $0$ otherwise.
Clearly $\|\tilde x\|_1 \le M$ and similarly for $\tilde x'$ yet we have that
$d(z, z') = 2M$. (This assumes that the number of entries, $n$, is even and
greater than one, otherwise one may scale the nonzero entries of $\tilde x'$,
since $\|\tilde x'\|_1 = M(1 - 2/n)$ when $n \ge 2$ is odd, to match the
bound.) 

\paragraph{Value function.} We assume that the market acts like a two-asset
fee-free constant function market maker~\cite{angeris2022}, which we define as
a nondecreasing concave function $G: \reals_+ \to \reals_+$ with $G(0) = 0$.
(This function $G$ is sometimes called the forward exchange
function~\cite[\S4.1]{angeris2022}.) We assume the market begins at some state
$0$ and, assuming some sequence of trades $t_1, \dots, t_k \in \reals_+$ were
made, then trade $k+1$, with some input amount $t_{k+1}$, pays out
\begin{equation}\label{eq:trades}
    G(t_1 + \dots + t_k + t_{k+1}) - G(t_1 + \dots + t_k).
\end{equation}
In this case, the payoff for the validator is as follows. Given some
actions $z = (\tilde x, \tilde y) \in B$, let $k$ be the index such that $\tilde y_k = \delta$,
then the payoff for the validator is
\[
    f(z) = G(\tilde x_1 + \dots + \tilde x_{k-1} + \delta) - G(\tilde x_1 + \dots + \tilde x_{k-1}).
\]

\paragraph{Smoothness.} Since $G$ is concave, we have that 
\[
    G'(t) \le G'(0),
\]
for $t \ge 0$, so we have, for $t' \ge 0$,
\[
    |G(t) - G(t')| \le G'(0)|t - t'|.
\]
This implies that, given two actions $z, z' \in B$ and setting
\[
    t = \sum_{i=1}^k \tilde x_i, \quad t' = \sum_{i=1}^{k'} \tilde x'_i,
\]
where $k$ is the index for which $\tilde y_k = \delta$, \ie, when the validator
makes their trade in action $z$, and $k'$ is similar for action $z'$, we have
\begin{multline}\label{eq:f-bound}
    |f(z) - f(z')| = |G(t + \delta) - G(t) - (G(t' + \delta) - G(t'))| \\
    \le |G(t + \delta) - G(t' + \delta)| + |G(t) - G(t')|  \le 2G'(0)|t -t'|.
\end{multline}
Using the triangle inequality, we have,
\[
    |t - t'| \le \sum_{i=1}^k |\tilde x_i - \tilde x_i'| + \sum_{i=k+1}^{k'} |\tilde x_i'| \le d(z, z').
\]
Here we assumed $k \le k'$ but swapping $z$ for $z'$ gives the complete claim.
We then get the final result using~\eqref{eq:f-bound} and the above:
\[
    |f(z) - f(z')| \le 2G'(0)d(z, z').
\]

\paragraph{Cost of MEV.} Using the bound~\eqref{eq:smooth-upper-bound} and the
`diameter' bound~\eqref{eq:diam-bound}, we then get that the cost of MEV is at
most
\[
    C_s(f) \le 8G'(0)M,
\]
where we used the fact that $\inf_{x\in B} C(f, x) = 0$ by setting all traders'
trades to be zero. This bound matches our intuitions: as the size of others'
trades grows larger ($M$ is large) the cost of MEV similarly becomes larger,
roughly linearly in the total amount of volume traded. (Roughly speaking, we
expect that the larger the total volume of others' trades, the more a validator
would be willing to pay to front-run.) This bound is only reasonable when the
validator's trade $\delta$ is at least on the order of the total trade volume
$M$; otherwise, the bound is very loose when $\delta$ is very small.

\subsubsection{Sandwiching}
Another example of MEV in decentralized exchanges is the notion of
\emph{sandwich attacks}, \ie, an adversary (such as the validator introduced
above) inserts a buy transaction before users purchase an asset and a sell
transaction after. Note that, intuitively, the amount of profit generated is
always nonnegative since the adversary always `purchases' an asset at a price
no higher than she sells it at.

\paragraph{Action space.} The set of possible actions in this setting takes the form
\[
    A = ( \reals_+ \times \{0\}) \cup ( \{0\} \times \reals )
\]
Similar to the previous example, a tuple $(x, y) \in A$ is an action where $x$
is a trade made by a trader, which is nonzero only when $y$ (the trade made by
the adversary) is zero. Alternatively, we can have either $y > 0$ or $y < 0$ as
the possible choices for the sandwicher's trade, which can be interpreted as
buy and sell orders in the decentralized exchange, respectively.

\paragraph{Sandwiching mechanics.} In our model, we will assume that the
sandwicher observes the trades to be included in the block, which we call $x_1,
\dots, x_k \in \reals_+$. The sandwicher then submits a buy order of
$\delta \in \reals_+$, ideally inserted before all other trades happen, and
submits a sell order of $\gamma \in \reals_+$, ideally inserted after all other
trades happen. We assume the same model of the market as the previous example,
defined by some function $G: \reals \to \reals$, which is strictly increasing and
concave, with $G(0) = 0$. (We allow negative trades in this case to indicate an
amount which must be tendered instead of received.) Assuming the sandwicher
purchases some amount of an asset by tendering $\delta$ and then sells back no
more than the amount $G(\delta)$ received, the trades must satisfy the
inequality
\begin{equation}\label{eq:sandwich}
    G(t + \delta) - G(t + \delta - \gamma) \le G(\delta) - G(0) = G(\delta),
\end{equation}
where $t = x_1 + x_2 + \dots + x_k$ is the total of all other trades. Note that
this uses the definition~\eqref{eq:trades} of the market $G$, and assumes that
the trade $\delta$ is made in the first position and $\gamma$ is made in the
last. (It is a useful exercise to show that, due to the concavity of $G$, and
the fact that all trades are nonnegative, this is the largest possible $\delta$
and $\gamma$ that can be tendered, or received, respectively.) The left hand
side can be interpreted as the amount of $\gamma$ received for tendering
$G(\delta)$ back to the market. The inequality simply says that this amount,
interpreted as the amount tendered, is bounded by the amount of the asset
received by the initial trade of $\delta$ (\ie, the sandwicher does not tend
more to the market than the amount that was received from the initial trade
of $\delta$).

\paragraph{Constrained action space.} We will define
the set of possible actions $B$ in the following way:
\begin{multline*}
    B = \{((x_1, y_1), \dots, (x_n, y_n)) \in A^n \mid \|x\|_1 \leq M ~\text{and exactly one} \\
    \text{pair $i,j$ with $y_i > 0$ and $y_j < 0$ which must satisfy~\eqref{eq:sandwich}} \}.
\end{multline*}
In other words, $B$ is defined as the set of trades for which the total trade
volume, excluding the sandwicher's trades, is no more than $M$, and the only
two trades by the sandwicher, given by, say, $y_i > 0$ and $y_j < 0$,
satisfy~\eqref{eq:sandwich} in the following way:
\begin{equation}\label{eq:sandwich-2}
    G(t + y_i) - G(t + y_i + y_j) \le G(y_i),
\end{equation}
where $t = x_1 + x_2 + \dots + x_n$.

\paragraph{Value function.} Overloading notation in the same way as the
previous example, we may write an action $z = (\tilde x, \tilde y) \in B$. We
define the value function, or payoff, for this action (from the perspective of
the sandwicher) to be the amount received of the asset minus the amount
tendered originally, which is simply
\begin{equation}\label{eq:sand-payoff}
    f(z) = -(\tilde y_j + \tilde y_i).
\end{equation}
(The negative sign comes from the fact that negative values denote amounts
received from the market.)

\paragraph{Metric.} Finally, we set the `metric' $d$ for $z = (\tilde{x},
\tilde{y}) \in B$ and $z' = (\tilde x', \tilde y') \in B$, where we have
overloaded notation in the same way as the previous example, to
\[
    d(z,z') = \max\{\|\tilde{x}\|_1, \|\tilde{x}'\|_1\}.
\]
In other words, this `metric' simply measures the largest volume of trades from
the two possible actions. By definition of $B$ we have that the `diameter'
is no more than
\[
    d(z, z') \le M.
\]

\paragraph{Smoothness.} Using the definition of $B$, the smoothness is very
easy to prove. Let $z = (\tilde x, \tilde y) \in B$ be an action. Note that, since
$G$ is increasing, then, using~\eqref{eq:sandwich-2} we have
\[
    G(t + \tilde y_i + \tilde y_j) \ge G(t + \tilde y_i) - G(t) \ge 0,
\]
where $t = \tilde x_1 + \dots + \tilde x_n$.
This implies that
\[
    t + \tilde y_i + \tilde y_j \ge 0,
\]
or that the payoff for any action is bounded from above by
\[
    f(z) = -(\tilde y_i + \tilde y_j) \le t \le \|\tilde x\|_1.
\]
Note that this inequality is rather intuitive: it states that
the maximum amount of MEV available from sandwiching cannot be
more than the total volume of trades happening on the market.
Since, from before, the payoff is always nonnegative, we have
that $f(z), f(z') \ge 0$ so,
\[
    |f(z) - f(z')| \le \max\{f(z), f(z')\} \le \max\{\|\tilde x\|_1, \|\tilde x'\|_1\} = d(z, z'),
\]
which means that $f$ is 1-smooth with respect to the `metric' $d$.

\paragraph{Cost of MEV.} Using~\eqref{eq:smooth-upper-bound},
and the fact that the `diameter' of $B$ is no more than $M$,
we find that the worst-case cost of MEV is
\[
    C_s(f) \le 2M.
\]
This immediately strengthens known bounds for $C_s(f)$
from~\cite{kulkarni2022towards}, which had $C_s(f) = O(\log n)$ versus $C_s(f)
= O(1)$, assuming the market volume is constant, as we have here. Again, this
is reasonably intuitive: the maximum cost of MEV cannot be more than a constant
factor away from the total market volume. In this particular case, using the
fact that $f$ is nonnegative, we may actually strengthen this result to
\[
    C_s(f) \le M,
\]
by directly using the definition of the cost of MEV.

\section{Permutations on graphs}\label{sec:permutations}
In this section, we will show the connections between the function $f(\pi(x))$
(where the permutation $\pi\in S$ is free while the transaction $x \in A^n$ is
fixed) and a graph whose vertices are indexed by the set of permutations $S$,
with edges encoding particular structure that we define next. We will show a
number of basic results which generalize some of those known in the literature.

\subsection{Definitions and examples}\label{sec:defininggraphs}
We will define a \emph{permutation graph}
as a connected, undirected graph $G=(V, E)$ where the vertices $V$ are indexed
by the set of permutations $S$ and the edge set $E$ is composed of subsets of
$V$ of cardinality $2$. While the vertices of the permutation graph are fixed
in this definition, the edges are not and we are free to choose $E$ in any way
we like, so long as the graph $G$ is connected. (We will make use of this
freedom later when analyzing $f$.) For convenience, we will overload notation
slightly by defining:
\[
    f_i = f(\pi_i(x)), \qquad i=1, \dots, n!,
\]
for fixed $x \in A^n$ and where $\pi_i$ ranges over all permutations. Unless
stated otherwise, we will assume that $x$ is fixed and implicitly included in
the definition of the $f_i$ for the remainder of the paper. We will also
identify $f \in \reals^{n!}$ as an $n!$-vector such that quantities like the
$\ell_2$ norm, $\|f\|_2$, have the obvious interpretations.

\paragraph{A basic bound.} While this definition appears trivial, there are already a number
of useful consequences. For example, for any path $P \subseteq E$ from node $i$
to node $j$,
\[
    |f_i - f_j| \le \sum_{\{k, \ell\} \in P} |f_k - f_\ell|,
\]
which follows from the triangle inequality. Taking the maximum over $i$ and $j$
and a shortest path for each pair, we have that
\begin{equation}\label{eq:l1-bound}
    \max_{i, j \in V} |f_i - f_j| \le \diam(G)\left(\max_{\{k, \ell\} \in E} |f_k - f_\ell|\right),
\end{equation}
where $\diam(G)$ is the diameter of the graph; \ie, the largest distance
between any two nodes of the graph, where the distance is defined as the length
$|P|$ of the shortest path $P$ between the nodes. (This is always finite as we
assumed the graph $G$ is connected.) Since
\[
    C(f) \le \max_{i, j \in V} |f_i - f_j|,
\]
this also immediately implies a bound on the cost of MEV.
Inequality~\eqref{eq:l1-bound} already gives us an idea for why a `good' choice
of $E$ is essential: in order to get a tight bound, we would like to choose $E$
to be `simple enough' to get a reasonable bound on the maximum difference
between adjacent permutations, while making it connected enough so the diameter
of the graph, $\diam(G)$, is small. We outline a few important examples below.

\paragraph{Complete graph.} Perhaps the simplest example of an edge set $E$ is
the set of all possible edges, in which case the graph $G$ is the complete
graph of $n!$ vertices. In this case bound~\eqref{eq:l1-bound} is trivial as
$\diam(G) = 1$ (since any node is exactly one hop away from any other node)
while
\[
    \max_{\{k, \ell\} \in E} |f_k - f_\ell| = \max_{k, \ell \in V} |f_k - f_\ell|,
\]
which is what we wanted to bound in the first place. In almost all cases, the
complete graph almost universally yields weak and/or trivial bounds for many
quantities, so it is only useful as an example of the definitions.

\paragraph{Transposition graph.}
A second possible choice of edge set $E$ is for every pair of nodes $i$ and
$j$, we add an edge $\{i, j\}$ to $E$ if $\pi_i$ and $\pi_j$ differ by exactly
one transposition. (In other words, if there are two elements in $\pi_i$ which,
when swapped, results in $\pi_j$.) We will call this graph the
\emph{transposition graph}. This graph is connected and has diameter $\diam(G)
= n-1$ since any permutation can be written as $n-1$ transpositions
(\cf,~\cite{conradSign}). If we have any bound on the maximum difference of $f$
in transpositions, say $|f_i - f_j| \le C$ if $\{i,j\} \in E$, then,
applying~\eqref{eq:l1-bound}, we get:
\begin{equation}\label{eq:perm-graph}
    \max_{i,j \in V} |f_i - f_j| \le (n-1)C.
\end{equation}
This bound strengthens the result
of~\cite{chitraDifferentialPrivacyConstant2022}, which considers the case where
the trades for a specific CFMM can be reordered to maximize some payoff, by a
factor of $\log(n)$. This graph is also known as the Cayley graph of the
symmetric group generated by the transpositions, and, as one might expect, has
deep connections to group theory (\cf,~\cite{diaconis1988group}). We point out
that the graph is $\binom{n}{2}$-regular (as there are that many possible
transpositions) and bipartite for future reference.


\subsection{Smoothness over a graph}

One simple way of thinking about the cost of MEV defined in~\eqref{eq:main} is
as a measure of `smoothness' of the function $f$ over the set of permutations
$S$, for a fixed set of transactions $x$. Given the graph $G = (V, E)$, we may
define a reasonable notion of `smoothness', as measured over the graph,
by setting
\[
    C_G(f) = \sum_{\{i,j\} \in E} (f_i - f_j)^2.
\]
The lower the value of $C_G(f)$ the more smooth a function $f$ is over all
`adjacent' permutations. Note that, for any $f$ this function is nonnegative
since it is the sum of nonnegative terms. We will see that $C_G$ and $C$ are
related, but $C_G$ can sometimes provide a simpler, and more intuitive, description
of `smoothness' over permutations.

\paragraph{Properties.} Like the original definition of $C$, the
value of $C_G$ is unchanged by offsets,
\ie,
\[
    C_G(f + \alpha \ones) = C_G(f),
\]
for $\alpha \in \reals$, where $\ones \in \reals^{n!}$ is the all-ones vector.
Additionally, $C_G$ is homogeneous of degree two in that
\[
    C_G(\alpha f) = \alpha^2 C_G(f),
\]
for any $\alpha \in \reals$.

\paragraph{Partial equivalence.} The suggestive notation $C_G$ here is in part
justified by the following simple fact: for a fixed set of transactions $x \in
A^n$, a function $f$ has $C_G(f) = 0$ if, and only if, $C(f) = C(f, x) = 0$.
(While the notation here is purposefully overloaded, the
interpretations of each statement should be clear from context and the fact
that $x$ is fixed.) To see the forward implication, note that $C_G(f) = 0$
implies that
\begin{equation}\label{eq:constantPayoff}
    f_i = f_j, \quad \text{for every} ~ \{i, j\} \in E.
\end{equation}
Since the graph is connected, then $f_i = f_1$ for every $i=1, \dots, n!$ and
therefore $f$ is a constant vector. Written another way, $f(\pi_i(x)) = f(x)$
for every permutation $\pi_i$, which exactly when the function is fair for
transactions $x$, \ie, when $C(f) = 0$. The reverse implication is immediate
from the definitions. (We will show a tighter connection between $C$ and $C_G$
later in this section.)

\subsubsection{The Laplacian of a graph}
The function $C_G(f)$ can be written as the following homogeneous quadratic:
\[
    C_G(f) = f^TLf,
\]
where the matrix $L \in \reals^{n! \times n!}$ is defined as
\[
    L_{ij} = \begin{cases}
        -1 & \{i, j\} \in E\\
        d_i & i = j,
    \end{cases}
\]
and $d_i$ is defined as the degree of node $i$; \ie, it is the number of
vertices adjacent to node $i$. This matrix $L$ is called the \emph{Laplacian}
and encodes a number of important properties of the graph $G$.

\paragraph{Eigenvectors and spectra.}
The matrix $L$ is symmetric since $\{i, j\} = \{j, i\}$. A basic fact from
linear algebra is that any symmetric matrix has an eigenvalue
decomposition~\cite{axler1997linear}:
\[
    L =U \Sigma U^T,
\]
where $\Sigma \in \reals^{n!\times n!}$ is a diagonal matrix, while $U \in
\reals^{n! \times n!}$ is an orthogonal matrix with $U^T U = UU^T = I$. For the
remainder of the paper, we will write $u_i \in \reals^{n!}$ for the $i$th
column of the matrix $U$ and $\lambda_i = \Sigma_{ii}$. Since we may rearrange
the columns of $U$ and $\Sigma$ as needed, we will also assume that $\lambda_1
\le \lambda_2 \le \dots \le \lambda_{n!}$. (This list of $\lambda_i$ is called
the \emph{spectrum} of the graph $G$.) Using the definition of $L$ and the
$u_i$ we have that
\[
Lu_i = \lambda_i u_i,
\]
(which is why the above decomposition is called an eigenvalue decomposition)
and therefore that
\begin{equation}\label{eq:eigvals}
    \lambda_i = u_i^TL u_i = C_G(u_i) \ge 0.
\end{equation}
This implies that the $\lambda_i$ are nonnegative.
We can show that $\lambda_1 = 0$ since, using the definition of $L$,
\[
    L\ones = 0\ones = 0,
\]
where $\ones$ is the all-ones vector, so we may assume that $u_1 = \ones/\sqrt{n!}$.
It is not hard to show that, since $G$ is connected, $\lambda_2 > 0$.

\paragraph{Bounds on the largest eigenvalue.} If the graph is a $k$-regular
graph, such as the transposition graph with $k=\binom{n}{2}$, then we have the
following bound on the largest eigenvalue:
\[
    \lambda_{n!} \le 2k.
\]
To see this, let $\tilde u$ correspond to any eigenvector with eigenvalue
$\tilde \lambda\ge 0$, and let $i \in \argmax_j |\tilde u_j|$ be its largest index, then:
\[
    \tilde \lambda |\tilde u_i| = |(L\tilde u)_i| = \Bigg|k\tilde u_i - \sum_{\{i, j\} \in E} \tilde u_j\Bigg| \le k|\tilde u_i| + \sum_{\{i, j\} \in E} |\tilde u_j| \le k|\tilde u_i| + k|\tilde u_i| = 2k|\tilde u_i|.
\]
The only major step is in the second inequality, which uses both the fact that
the graph is $k$-regular (\ie, every node has exactly $k$-neighbors) and the
fact that $|\tilde u_i|$ is the largest entry in $\tilde u$. Since $\tilde u
\ne 0$ then $\tilde u_i \ne 0$ by definition and canceling both sides we have
\[
    \tilde \lambda \le 2k,
\]
as required. (Using $\tilde \lambda = \lambda_{n!}$ and similarly for $\tilde
u$ gives the final result.) We will make use of this result, applied to the
special case of the transposition graph, later in this section.

\subsubsection{Fourier transform over graphs}
We are now ready to define the \emph{Fourier transform} of the function $f$
over the graph $G$, given simply by the vector $\hat f \in \reals^{n!}$,
defined
\[
    \hat f = U^Tf = \begin{bmatrix}
        u_1^Tf\\
        u_2^Tf\\
        \vdots\\
        u_{n!}^Tf
    \end{bmatrix},
\]
where $u_i$ is the $i$th eigenvector of the Laplacian $L$. Since $U U^T = I$ we
can similarly write
\begin{equation}\label{eq:decomp}
    f = U(U^Tf) = U\hat f = \sum_{i=1}^{n!} \hat f_i u_i.
\end{equation}
In this notation, it is clear that $u_i$ act as a \emph{Fourier basis} and each of the $\hat{f}_i$ act as \emph{Fourier coefficients} that weight the basis vectors. Indeed, from ~\eqref{eq:eigvals} and their definition, the $u_i$ have a reasonably
simple interpretation: the $u_i$ form an orthogonal basis with the property
that, as $i$ becomes larger, the $u_i$ become `less smooth' over $G$, or, in
signal processing parlance, they represent `higher frequencies' of the graph.
We may then decompose $f$ into its constituent `frequencies', where $\hat f_i$
gives the contribution of each frequency $v_i$ to the overall function $f$.
(Note that $u_1$ represents the lowest possible frequency: a constant.)

In this sense, the decomposition of~\eqref{eq:decomp} is very similar (and, in
fact, in certain limits, for certain graphs, equivalent~\cite{singer2006}) to
the traditional Fourier transform for signals. We have many of the usual
identities, such as Parseval's theorem, hold:
\begin{equation}\label{eq:parseval}
    \|\hat f\|_2^2 = \hat f^T\hat f = \hat f^T(U^T U)\hat f = (U\hat f)^T(U\hat f) = \|f\|_2^2,
\end{equation}
since $U^T U = I$. We will also see next that the Fourier transform of $f$ over
$G$ has deep connections with the definition of fairness presented in~\eqref{eq:main}.
As an important aside: the transform $\hat f$ depends heavily on the structure
of the graph $G$ and may look very different for different graphs (even though
the function $f$ may be unchanged).

\paragraph{Largest eigenvector.} For the case of the 
transposition graph, it is possible to explicitly construct its largest eigenvector.
To do this, first note that the transposition graph is $\binom{n}{2}$-regular.
This means that its largest eigenvalue, $\lambda_{n!}$, is bounded by
\[
    \lambda_{n!} \le 2\binom{n}{2} = n(n-1),
\]
from the previous discussion. Now, since, for any $v \in \reals^{n!}$,
\[
    C_G(v) = v^TLv \le \lambda_{n!}\|v\|_2^2,
\]
we will exhibit a vector $v$ with cost $C_G(v) = n(n-1)\|v\|_2^2$, which will
show that $\lambda_{n!} = n(n-1)$ and that $v$ is a corresponding eigenvector.
We will define $v$ in the following way
\[
    v_i = \begin{cases}
        +1 & \text{if $\pi_i$ is an even permutation}\\
        -1 & \text{if $\pi_i$ is an odd permutation},
    \end{cases}
\]
where we have assumed that the $i$th node corresponds to permutation $\pi_i \in
S$, for $i=1, \dots, n!$. (For much more on the definition and existence of
even/odd permutations, see~\cite{conradSign}.) Now, note that
\[
    C_G(v) = \sum_{\{i,j\} \in E} (v_i - v_j)^2.
\]
But, by definition of the transposition graph, two nodes $i$, $j$, are
neighbors if, and only if, they differ by a transposition. So, if $\pi_i$ is
odd, then $\pi_j$ must be even or vice versa, hence
\[
    C_G(v) = |E|(2)^2 = 4|E|,
\]
where $|E|$ is the number of edges. Since the transposition graph is
$\binom{n}{2}$-regular, we know that
\[
    \frac{n!\binom{n}{2}}{2} = |E|,
\]
and we also know that $\|v\|_2^2 = n!$, so, we get
\[
    C_G(v) = 4|E| = 2n!\binom{n}{2} = n(n-1)\|v\|_2^2,
\]
as required. Normalizing this function gives an eigenvector corresponding
to the largest eigenvalue,
\[
    u_{n!} = \frac{1}{\sqrt{n!}}v.
\]
(In fact, the largest eigenvalue has multiplicity one, so this eigenvector is
unique, but we do not prove this here.)

\paragraph{Graph Fourier coherence.} Having defined the Fourier transform of a 
graph, we now define a quantity that will be crucial in giving bounds on the 
cost of MEV, following \cite{perraudin2018global}. For every set of orthogonal 
eigenvectors $u_j$ corresponding to a graph $G$, we define the \emph{graph Fourier 
coherence} as:
\begin{equation}\label{eq:coherence}
    \mu = \max_i \|u_i\|_\infty,
\end{equation}
where $\|y\|_\infty = \max_i |y_i|$. This quantity measures how `concentrated'
the Fourier basis is when compared to the standard basis. Using basic estimates
for the $\ell_\infty$ norm in the $\ell_2$ norm since $\|u_i\|_2 = 1$:
\[
    \frac{1}{\sqrt{n!}} = \frac{1}{\sqrt{n!}}\|u_i\|_2 \le \|u_i\|_\infty \le \|u_i\|_2 = 1,
\]
so $1/\sqrt{n!} \le \mu \le 1$. The upper bound occurs precisely when
there is a vector $u_i$ that is exactly aligned with one of the standard basis
vectors, while the minimum value indicates that the Fourier basis is maximally
`incoherent' with the standard basis. As we will see, bounds on this quantity
control our ability to localize a function $f$ and its Fourier coefficients
$\hat{f}$ simultaneously. Intuitively, if $\mu$ is small, then the `size' of
$f$ and $\hat{f}$ must be traded off against each other in a sense that we will
make precise. Note that, because $\mu$ is a graph-dependent quantity, the way
that the edge set is constructed will affect this quantity, and hence our
ability to bound the cost of MEV. Note that the coherence~\eqref{eq:coherence} is
a basis dependent phenomenon: if the eigenvalues are degenerate (\ie, if there
are repeated eigenvalues) then the coherence need not be unique and would
depend on the specific eigenvectors chosen in the decomposition.

\paragraph{Coherence for connected graphs.} In the case where a graph is
connected and has more than one vertex, we can get a slightly better upper
bound on the coherence of any eigenvector $v$, given by
\[
    \|v\|_\infty \le \sqrt{1 - \frac{1}{n!}}.
\]
In particular, the first eigenvector is known (the all-ones vector) and has
\[
    \frac{\|\ones\|_\infty}{\sqrt{n!}} = \frac{1}{\sqrt{n!}},
\]
which is exactly the lower bound on $\mu$ given above. Since the graph is
connected, we must have that every other eigenvector $v$ satisfies $\ones^Tv =
0$, from the previous discussion. For any nonzero vector, there must be at
least one component that is strictly positive and one that is strictly
negative; we write $v^+$ for the positive entries of $v$ and $v^-$ for the
negative entries and assume, without loss of generality, that the maximum
absolute value entry is in the positive part. (This is without loss of
generality since, if $v$ is an eigenvector, then $-v$ is as well.)
With this, we have
\[
    \|v\|_\infty \le \|v^-\|_1 \le \|v^-\|_2\sqrt{n!-1} = \sqrt{1-\|v^+\|_2^2}\sqrt{n!-1} \le \sqrt{1-\|v\|_\infty^2}\sqrt{n!-1}.
\]
Rearranging and solving for $\|v\|_\infty$ gives the inequality above.
The first inequality follows from the fact that the sums of the positive
and negative entries must be equal, so the largest entry (which is positive)
must be at least the sum of all of the negative entries. The second inequality
follows from the fact that there are at most $n!-1$ nonzero entries in $v^-$,
while the third follows since $\|v\|_2 = 1$. Finally, the last follows
from the fact that $\|v^+\|_\infty \le \|v^+\|_2$ and that the largest
entry is positive.

Note that this inequality is tight in that it can be saturated by choosing a
vector with $v_1 = \sqrt{1-1/n!}$ and $v_i = -1/\sqrt{n!(n!-1)}$ for $i=2, \dots, n!$.

\paragraph{Coherence for the complete graph.} We note that coherence is a
phenomenon that depends on the specific Fourier basis that is chosen. (In
particular, some graphs may have degenerate eigenvalues, in which case there
are many possible eigendecompositions.) In the case of the complete graph, we
have that there exists an eigendecomposition that saturates the above bound.
This follows easily from the fact that any vector orthogonal to the all-ones
vector is an eigenvector of the complete graph Laplacian. Choosing the vector
above that saturates this inequality then gives us that, indeed, the coherence
in any basis containing any such vector saturates the inequality for the
complete graph.

\paragraph{Numerical estimates.} We numerically compute $\mu$ for the various
kinds of permutation graphs that are considered in~\S\ref{sec:defininggraphs}.
The results for complete graphs and transposition graphs are provided in
table~\ref{table:fourier}. Further numerical estimates for other kinds of
graphs, including Erd\H{o}s-R\'enyi graphs can be found in
\cite{perraudin2018global}. It can be seen that as the number of transactions
increases, the coherence decays the fastest for transposition graphs; this
shows that the choice of the underlying graph can meaningfully affect the
bounds on the cost of MEV, as we will see next. Code to reproduce these
quantities is provided at
\begin{center}
    \texttt{https://github.com/bcc-research/spectrum-calculations}
\end{center}
Note that, in general, one might receive different numbers than those given in
the table, as the transposition graph has degenerate eigenvalues. (The complete
graph, on the other hand, is a global bound from the previous discussion.) The
Fourier coefficients therefore also depend on the exact eigenvectors chosen
which is why the graph coherence varies.

\begin{table}
    \centering
    \begin{tabular}{ccc}
        \toprule
        \textbf{Number of transactions, $n$} & \textbf{Transposition graph} & \textbf{Complete graph} \\
        \midrule
        1 & 1.000 & 1.000 \\
        2 & 0.707 & 0.707 \\
        3 & 0.816 & 0.913 \\
        4 & 0.612 & 0.978 \\
        5 & 0.548 & 0.995 \\
        6 & 0.323 & 0.999 \\
        7 & 0.477 & 1.000 \\
        \bottomrule
    \end{tabular}
    \caption{Graph Fourier coherence, $\mu$, for transposition and complete graphs.}
    \label{table:fourier}
\end{table}

\subsection{Fairness and the Fourier transform}
For the next section, we will assume that $\max_i f_i \ge \max_j -f_j$ ; in
other words, the highest possible payoff is larger than the negative of the
largest possible loss, which, using the translation invariance of $C$, may
generally be assumed. 

\paragraph{Cost of MEV.} Given the above condition on $f$, note that $\max_i
f_i = \|f\|_\infty$. Additionally, we may write the expectation of $f$
uniformly over the possible permutations as $\ones^Tf/n!$. This means that we
may rewrite the cost of MEV in~\eqref{eq:main} as (dropping the $x \in A^n$ as
it is fixed):
\[
    C(f) = \|f\|_\infty - \frac{1}{n!}\ones^Tf.
\]
Since $u_1 = \ones/\sqrt{n!}$ and $f$ can be decomposed as given in~\eqref{eq:decomp},
then
\[
    \frac{1}{n!}\ones^Tf =\frac{1}{n!} \sum_{i=1}^{n!} \ones^T u_i \hat f_i,
\]
but, since $u_1^Tu_i = \ones^Tu_i = 0$ unless $i=1$ we then have
\[
    \frac{1}{n!}\ones^Tf = \frac{1}{\sqrt{n!}} \hat f_1,
\]
so we can write
\[
    C(f) = \|f\|_\infty - \frac{1}{\sqrt{n!}}\hat f_1.
\]

\paragraph{Lower bounds.} This rewriting immediately provides a simple
lower bound only in terms of $\hat f$, since $\sqrt{n!}\|f\|_\infty \ge \|f\|_2 = \|\hat f\|_2$,
where the latter equality follows from~\eqref{eq:parseval}, so
\begin{equation}\label{eq:bound-upper}
    C(f) \ge \frac{1}{\sqrt{n!}}\left(\|\hat f\|_2 - \hat f_1\right).
\end{equation}
We can view the term $\|\hat f\|_2 - \hat f_1$ roughly as asking `how much more
do the higher frequencies of $f$ contribute over $\hat f_1$'? Alternatively:
how much of the mass of $\hat f$ is contained in all of the nonconstant
frequencies? This also gives an interesting observation. If $f$ is a perfectly
fair function, we know that $C(f) = 0$. Using~\eqref{eq:bound-upper}, we get
$\|\hat f\|_2 \le \hat f_1$, or, after squaring both sides and subtracting,
\[
    \hat f_2^2 + \hat f_3^2 + \dots + \hat f_{n!}^2 \le 0,
\]
which immediately implies that $\hat f_1$ is the only possible nonzero
coefficient. This corresponds neatly to the `usual' notion that a constant
function has a Fourier transform whose only support is at the lowest possible
frequency. (We show the `opposite', that worst-cost functions have
nonzero coefficients at high frequencies, later.)

\paragraph{Upper bounds.} A basic upper bound for this problem follows from
the fact that
\[
    \|f\|_\infty = \left\|\sum_{i=1}^{n!} \hat f_i u_i\right\|_\infty \le \sum_{i=1}^{n!} |\hat f_i| \|u_i\|_\infty \le \|\hat f\|_1\left(\max_i \|u_i\|_\infty\right) = \|\hat f\|_1 \mu,
\]
which gives
\[
    C(f) \le \mu \|\hat f\|_1 - \frac{1}{\sqrt{n!}}\hat f_1.
\]

\paragraph{Worst-cost functions on the transposition graph.} In a very specific
sense, the `worst-cost' functions defined in the previous sections are some of
the `spikiest' functions over the graph. For example, in the special case of
the transposition graph, the largest eigenvector is always nonzero at every
node of the graph, which immediately shows that
\[
    \hat e_i = e_i^Tu_{n!} = \pm 1/\sqrt{n!},
\]
so worst-cost functions always have nonzero components at the highest
possible frequency. Much in the same way that the (usual) Fourier transform of
very `spiky' functions has high frequency components, we see that the graph
Fourier transform of a worst-cost function (namely, the basis vectors)
similarly has high frequency components when considering the transposition
graph. 

\subsubsection{Relationship between $C$ and $C_G$}
Unsurprisingly, there is a relationship between the cost of MEV $C$ and
the measure of `smoothness' $C_G$, given by,
\[
    \frac{1}{\sqrt{\lambda_{n!} n!}} \sqrt{C_G(f)} \le C(f) \le \sqrt{\frac{C_G(f)}{\lambda_2}}.
\]
We show the two inequalities that relate these two functions in what follows.

\paragraph{Lower bound.} There is a lower bound in terms of
the largest eigenvalue given by
\[
    C(f) \ge \frac{1}{\sqrt{\lambda_{n!} n!}} \sqrt{C_G(f)}.
\]
To see that the original bound is true, note that, since $\lambda_{n!}$ is the
largest eigenvalue of $L$, we have
\[
    C_G(f) = f^TLf \le \lambda_{n!}\|f\|_2^2.
\]
Using~\eqref{eq:bound-upper} and the fact that $\|f\|_2 = \|\hat f\|_2$
we then have
\[
    C(f) \ge \frac{1}{\sqrt{n!}}\left(\sqrt{\frac{C_G(f)}{\lambda_{n!}}} - \hat f_1\right).
\]
Finally, using the fact that both $C$ and $C_G$ are translation invariant, set
$\alpha = \hat f_1/\sqrt{n!}$ and note that the first Fourier coefficient of
the function $f - \alpha \ones$ is zero to get
\[
    C(f) = C(f - \alpha\ones) \ge \frac{1}{\sqrt{n!}}\sqrt{\frac{C_G(f - \alpha\ones)}{\lambda_{n!}}} = \frac{1}{\sqrt{\lambda_{n!}n!}}\sqrt{C_G(f)}.
\]

\paragraph{Upper bound.} We also have an upper bound when $\ones^Tf \ge 0$
given by
\[
    C(f) \le \sqrt{\frac{C_G(f)}{\lambda_2}}.
\]
This is a refinement of the previous claim that $C_G(f) = 0$ implies that $C(f)
= 0$, since we know that $\lambda_2 > 0$ as the graph is connected. (Note that
this inequality fails when the graph is not connected as we would have
$\lambda_2 = 0$.) To see this, bound $C(f)$ by
\[
    C(f) = \max_i |f_i| - \ones^Tf/n! \le \max_i |f_i - \ones^Tf/n!| \le \|f - \alpha \ones\|_2,
\]
where we have defined $\alpha = \ones^Tf/n!$. We will overload notation
slightly and write $\hat f$ for the Fourier coefficients of  $f - \alpha\ones$.
The first coefficient of this function, $\hat f_1$, is zero since it is
orthogonal to the all-ones vector by definition of $\alpha$, $\ones^T(f -
\alpha\ones) = 0$. This means that
\[
    C_G(f - \alpha \ones) = \sum_{i=2}^{n!} \lambda_i \hat f_i^2 \ge \lambda_2 \sum_{i=2}^{n!} \hat f_i^2 = \lambda_2 \|\hat f\|_2^2,
\]
as the eigenvalues are in nondecreasing order, by assumption, and $\hat f_1 =
0$. Finally, using the fact that $\|\hat f\|_2 = \|f - \alpha\ones\|_2$ and our
previous inequality for $C$, we get
\[
    C(f) \le \sqrt{\frac{C_G(f-\alpha\ones)}{\lambda_2}} = \sqrt{\frac{C_G(f)}{\lambda_2}},
\]
where the last equality follows from the fact that $C_G$ is
translation-invariant. There is a (relatively weak) lower bound for $\lambda_2$
due to~\cite{mohar1991eigenvalues}, which is
\[
    \lambda_2 \ge \frac{4}{\diam(G)n!}.
\]
This would imply the following bound, which does not depend
on the eigenvalues of the Laplacian:
\[
    C(f) \le \frac{\sqrt{\diam(G)n!}}{2}\sqrt{C_G(f)}.
\]

\paragraph{Transposition graph.} In the special case of the transposition
graph, we know that it is $\binom{n}{2}$-regular and has diameter $n-1$, which
gives the following relationship between $C$ and $C_G$,
\[
    \frac{1}{\sqrt{n(n-1)n!}}\sqrt{C_G(f)} \le C(f) \le \frac{\sqrt{(n-1)n!}}{2}\sqrt{C_G(f)}.
\]
Note that these bounds are particularly loose and differ by large constants.
Since the transposition graph has a lot of structure that we do not exploit
here, it is very likely that the bounds provided can be significantly improved.

\section{Conclusion}
We constructed a simple framework that may be used to assess the economic
impact of MEV, which is excess value that a monopolistic validator
can extract from reordering transactions in a decentralized network. 
We defined the \emph{cost of MEV},
which is the difference between the worst-case and average case payoff to the
validator. Our definition allows us to compute a notion of `fairness' for a set
of transaction orderings that is a function of the economic value that can be
extracted from users. This is in contrast to definitions of `fair ordering' that
define `fairness' according to arrival times of transactions.

The key object in the definition is a payoff function: a real-valued function
that maps an ordered list of transactions to a utility, or payoff, realized by the
validator. We demonstrate that all payoff functions realized by users can be
written as linear combinations of the worst-cost payoff functions, which are
similar in spirit to liquidations---core mechanisms within decentralized
finance. We then make the distinction between `spiky' and `flat' payoff functions 
via simple bounds. We show that for functions that are smooth, so too is the
difference between their worst and average case behavior. 
Additionally, we
noted that payoff functions on permutations may also be seen analogously as
functions defined on graphs, where the vertices are the permutations. In this
case, we showed that the spectra of these functions provide non-trivial upper
and lower bounds on the cost of MEV, analyzed via the Fourier transform on
graphs.

\paragraph{Future work.} To better model new systems in production, it would be
interesting to consider `hierarchical' costs of MEV. A number of new MEV
auction designs such as SUAVE~\cite{MEV_Suave}, Jito~\cite{jitoJitoBlock},
Anoma~\cite{intentsanoma}, and Skip~\cite{skipXbuilderFirst} have auctions on a
per-domain or application level. A domain corresponds to a particular set of
contracts or code where reordering is allowed, while reordering between domains
is not controlled by any individual party (such as a monopolist sequencer or
validator)~\cite{obadia2021unity}. In such a world, one can imagine the set of
$n$ pending transactions as distributed into $k$ groups of transactions,
representing the $k$ distinct domains. One can compute a cost of MEV for each
group, but this naturally leads to a question: `what is the aggregate cost of
MEV over all domains?'. This open problem may have some analogues to the
current framework, which we leave for future work.

\section*{Acknowledgments}
The authors would like to thank Peteris Erins for suggestions, comments, and
edits on earlier drafts of this paper.

\bibliographystyle{alpha}
\bibliography{citations.bib}

\newcommand{\etalchar}[1]{$^{#1}$}
\begin{thebibliography}{MDFO22}

\bibitem[AAE{\etalchar{+}}22]{angeris2022}
Guillermo Angeris, Akshay Agrawal, Alex Evans, Tarun Chitra, and Stephen Boyd.
\newblock Constant {{Function Market Makers}}: {{Multi-asset Trades}} via
  {{Convex Optimization}}.
\newblock In {\em Handbook on {{Blockchain}}}, page~31. {Springer}, first
  edition, 2022.

\bibitem[Axl97]{axler1997linear}
Sheldon Axler.
\newblock {\em Linear algebra done right}.
\newblock Springer Science \& Business Media, 1997.

\bibitem[BCLL22]{bartoletti2022maximizing}
Massimo Bartoletti, James Hsin-yu Chiang, and Alberto Lluch~Lafuente.
\newblock Maximizing extractable value from automated market makers.
\newblock In {\em International Conference on Financial Cryptography and Data
  Security}, pages 3--19. Springer, 2022.

\bibitem[BO22]{bebel2022ferveo}
Joseph Bebel and Dev Ojha.
\newblock Ferveo: Threshold decryption for mempool privacy in bft networks.
\newblock {\em Cryptology ePrint Archive}, 2022.

\bibitem[CAE22]{chitraDifferentialPrivacyConstant2022}
Tarun Chitra, Guillermo Angeris, and Alex Evans.
\newblock Differential {{Privacy}} in {{Constant Function Market Makers}}.
\newblock {\em Financial Cryptography and Data Security 2022}, 2022.

\bibitem[CK22]{chitra2022improving}
Tarun Chitra and Kshitij Kulkarni.
\newblock Improving proof of stake economic security via mev redistribution.
\newblock In {\em Proceedings of the 2022 ACM CCS Workshop on Decentralized
  Finance and Security}, pages 1--7, 2022.

\bibitem[Cla84]{clark1984elements}
Allan Clark.
\newblock {\em Elements of abstract algebra}.
\newblock Courier Corporation, 1984.

\bibitem[Con]{conradSign}
Keith Conrad.
\newblock The sign of a permutation.
\newblock {\em Online}.

\bibitem[DGK{\etalchar{+}}20]{daian2020flash}
Philip Daian, Steven Goldfeder, Tyler Kell, Yunqi Li, Xueyuan Zhao, Iddo
  Bentov, Lorenz Breidenbach, and Ari Juels.
\newblock Flash boys 2.0: Frontrunning in decentralized exchanges, miner
  extractable value, and consensus instability.
\newblock In {\em 2020 IEEE Symposium on Security and Privacy (SP)}, pages
  910--927. IEEE, 2020.

\bibitem[Dia88]{diaconis1988group}
Persi Diaconis.
\newblock Group representations in probability and statistics.
\newblock {\em Lecture notes-monograph series}, 11:i--192, 1988.

\bibitem[Fla22]{MEV_Suave}
Flashbots.
\newblock The future of mev is suave: Flashbots, Nov 2022.

\bibitem[{Fla}23]{flashbots_mev_explore}
{Flashbots Team}.
\newblock Mev explore, 2023.

\bibitem[GSBV]{greeneoracle}
Jacob Greene, Ashar Shahid, Burak Benligiray, and Heikki V{\"a}nttinen.
\newblock Oracle extractable value (oev) through order flow auctions.

\bibitem[HPM23]{skipXbuilderFirst}
Sam Hart, Barry Plunkett, and Maghnus Mareneck.
\newblock x/builder: {T}he first {S}overeign {M}{E}{V} module for
  protocol-owned building --- ideas.skip.money.
\newblock
  \url{https://ideas.skip.money/t/x-builder-the-first-sovereign-mev-module-for-protocol-owned-building/57},
  2023.

\bibitem[JDE{\etalchar{+}}23]{johnson2023concave}
Nicholas~AG Johnson, Theo Diamandis, Alex Evans, Henry de~Valence, and
  Guillermo Angeris.
\newblock Concave pro-rata games.
\newblock {\em arXiv preprint arXiv:2302.02126}, 2023.

\bibitem[KDC22]{kulkarni2022towards}
Kshitij Kulkarni, Theo Diamandis, and Tarun Chitra.
\newblock Towards a theory of maximal extractable value i: Constant function
  market makers.
\newblock {\em arXiv preprint arXiv:2207.11835}, 2022.

\bibitem[KDL{\etalchar{+}}21]{kelkar2021themis}
Mahimna Kelkar, Soubhik Deb, Sishan Long, Ari Juels, and Sreeram Kannan.
\newblock Themis: Fast, strong order-fairness in byzantine consensus.
\newblock {\em Cryptology ePrint Archive}, 2021.

\bibitem[Lab23]{jitoJitoBlock}
Jito Labs.
\newblock {J}ito {B}lock {E}ngine {E}xpands {A}ccess to {A}ll {S}olana
  {M}{E}{V} {T}raders {J}ito {L}abs --- jito.wtf.
\newblock
  \url{https://www.jito.wtf/blog/jito-block-engine-expands-access-to-all-solana-mev-traders/},
  2023.

\bibitem[Mac23]{macpherson2023adversarial}
Andrew~W Macpherson.
\newblock Adversarial queues and trading on a cfmm.
\newblock {\em arXiv preprint arXiv:2302.01663}, 2023.

\bibitem[MDFO22]{mcmenamin2022fairtradex}
Conor McMenamin, Vanesa Daza, Matthias Fitzi, and Padraic O'Donoghue.
\newblock Fairtradex: A decentralised exchange preventing value extraction.
\newblock In {\em Proceedings of the 2022 ACM CCS Workshop on Decentralized
  Finance and Security}, pages 39--46, 2022.

\bibitem[Moh91]{mohar1991eigenvalues}
Bojan Mohar.
\newblock Eigenvalues, diameter, and mean distance in graphs.
\newblock {\em Graphs and combinatorics}, 7(1):53--64, 1991.

\bibitem[MS22]{malkhi2022maximal}
Dahlia Malkhi and Pawel Szalachowski.
\newblock Maximal extractable value (mev) protection on a dag.
\newblock {\em arXiv preprint arXiv:2208.00940}, 2022.

\bibitem[OSS{\etalchar{+}}21]{obadia2021unity}
Alexandre Obadia, Alejo Salles, Lakshman Sankar, Tarun Chitra, Vaibhav
  Chellani, and Philip Daian.
\newblock Unity is strength: A formalization of cross-domain maximal
  extractable value.
\newblock {\em arXiv preprint arXiv:2112.01472}, 2021.

\bibitem[PRSV18]{perraudin2018global}
Nathanael Perraudin, Benjamin Ricaud, David~I Shuman, and Pierre Vandergheynst.
\newblock Global and local uncertainty principles for signals on graphs.
\newblock {\em APSIPA Transactions on Signal and Information Processing}, 7,
  2018.

\bibitem[QZG22]{qin2022quantifying}
Kaihua Qin, Liyi Zhou, and Arthur Gervais.
\newblock Quantifying blockchain extractable value: How dark is the forest?
\newblock In {\em 2022 IEEE Symposium on Security and Privacy (SP)}, pages
  198--214. IEEE, 2022.

\bibitem[Res23]{resnick2023contingent}
Max Resnick.
\newblock Contingent fees in order flow auctions.
\newblock {\em arXiv preprint arXiv:2304.04981}, 2023.

\bibitem[RK23]{rondelet2023threshold}
Antoine Rondelet and Quintus Kilbourn.
\newblock Threshold encrypted mempools: Limitations and considerations.
\newblock {\em arXiv preprint arXiv:2307.10878}, 2023.

\bibitem[Rou16]{roughgarden2016}
Tim Roughgarden.
\newblock {\em Twenty Lectures on Algorithmic Game Theory}.
\newblock Cambridge University Press, 2016.

\bibitem[Sin06]{singer2006}
A.~Singer.
\newblock From graph to manifold laplacian: The convergence rate.
\newblock {\em Applied and Computational Harmonic Analysis}, 21(1):128--134,
  2006.
\newblock Special Issue: Diffusion Maps and Wavelets.

\bibitem[SY22]{intentsanoma}
Awa Sun~Yin.
\newblock An overview of anoma’s architecture.
\newblock Anoma, 2022.

\bibitem[XFP23]{xavier2023credible}
Matheus~Venturyne Xavier~Ferreira and David~C Parkes.
\newblock Credible decentralized exchange design via verifiable sequencing
  rules.
\newblock In {\em Proceedings of the 55th Annual ACM Symposium on Theory of
  Computing}, pages 723--736, 2023.

\bibitem[ZQT{\etalchar{+}}21]{zhou2021high}
Liyi Zhou, Kaihua Qin, Christof~Ferreira Torres, Duc~V Le, and Arthur Gervais.
\newblock High-frequency trading on decentralized on-chain exchanges.
\newblock In {\em 2021 IEEE Symposium on Security and Privacy (SP)}, pages
  428--445. IEEE, 2021.

\end{thebibliography}

\end{document}